\documentclass[leqno,10pt]{article}
\usepackage{amsmath}
\usepackage{amssymb,amscd}
\usepackage{latexsym}

\usepackage{color} 
 \newcommand{\ec}{\color{black}} 
\newtheorem{thm}{Theorem}[section]

\newtheorem{prop}[thm]{Proposition}
\newtheorem{cor}[thm]{Corollary}
\newtheorem{rem}[thm]{Remark}

\newcommand{\RE}{R \!\Join\! E}
\newcommand{\RI}{R\!\Join\! I}

\newcommand{\Hom} {\mbox{\rm Hom}}

\newcommand{\calF} {\mathcal F}

\title{The amalgamated duplication of a ring  \\  along a multiplicative-canonical ideal
}
\author{Marco D'Anna{\footnote{Partially supported by MIUR, under Grant PRIN 2005-011955}}
\and Marco Fontana{\footnote{Partially supported by MIUR, under Grant PRIN 2005-015278}}}

\begin{document}

\maketitle

\begin{abstract}

After recalling briefly the main properties of the amalgamated
duplication of a ring $R$ along an ideal  $I$, denoted by
$R\!\Join\!I$ \cite{d'a-f-1}, we restrict our attention to the
study of the properties of $R\!\Join\!I$,  when $I$ is a
multiplicative canonical ideal of $R$  \cite{hhp}.   In
particular, we study when every regular fractional ideal of
$R\!\Join\! I$ is divisorial.

\end{abstract}

\bigskip

%

\section{Introduction}

If $R$ is a commutative ring with unity and $E$ is an $R$-module,
\sl the idealization $R \!\ltimes\! E$ \rm  (also called \sl
trivial extension\rm ), introduced by Nagata in 1956 (cf. Nagata's
book \cite{n}, page 2), is a new ring where the module $E$ can be
viewed as an ideal such that its square is $(0)$. This
construction has been used in many contexts to produce examples of
rings satisfying preassigned conditions (see e.g. Huckaba's book
\cite{h}).
 In particular, in \cite[Theorem 7]{r} Reiten
proved that, if $R$ is a local Cohen Macaulay ring, then $R \!
\ltimes\! E$ is Gorenstein if and only if   $E$  is a canonical
module of $R$ (cf. also \cite[Theorem 5.6]{fgr}).

Fossum, in \cite{fs}, generalized the idealization defining \sl  a
commutative extension \rm of a ring $R$ by an $R$--module $E$ and
proved that, if $R$ is a local Cohen-Macaulay ring and if $E$ is a
canonical module of $R$, then any commutative extension $S$ of $R$
by $E$ is a Gorenstein ring \cite[Theorem]{fs}.

In this paper, we deal with some applications of a similar general
construction, introduced recently in \cite{d'a-f-1}, called \sl
the amalgamated duplication of a ring $R$ along an $R$--module
$E$, \rm which is an ideal in some overring of $R$, and denoted by
$R\!\Join\!E$. When $E^2 =0$, the new construction $R\!\Join\!E$
coincides with the Nagata's idealization $R\!\ltimes\! E$. In
general, however,  $\RE$  is not a commutative extension in the
sense of Fossum. One main difference of this construction, with
respect to the idealization (or with respect to any commutative
extension, in the sense of Fossum) is that the ring $\RE$ can be a
reduced ring (and it is always reduced if $R$ is a domain).

As it happens for the idealization, one interesting application of
this construction is the fact that it allows to produce rings
satisfying (or not satisfying) preassigned conditions. Moreover,
in many cases, the amalgamated duplication of a ring preserves the
property of being reduced (see \cite{d'a}, \cite{d'a-f-1}). Note
also that this new construction has been already applied for
studying questions concerning the diameter and girth of the
zero-divisor graph of a ring  (see \cite{m-y}).

M. D'Anna \cite{d'a} has studied this construction in case $E=I$
is a proper ideal of $R$, proving that, if $R$ is a local
Cohen-Macaulay ring with canonical module $\boldsymbol{\omega}_R$,
then $\RI$ is a Gorenstein ring if and only if $I \cong
\boldsymbol{\omega}_R$. \ec

Since in the one-dimensional local Cohen-Macaulay case the
Gorenstein rings are characterized by the property that  the
regular ideals are divisorial, it is natural to ask in a general
(non necessarily Noetherian) setting,  when $I$ is a
multiplicative canonical ideal of $R$, whether every regular
fractional ideal of the ring $\RI$ is divisorial.  Recall that the
notion of multiplicative canonical ideal was introduced in the
integral domain case by W. Heinzer, J. Huckaba and I. Papick
\cite{hhp}, and it can be easily extended  to  any commutative
ring: a regular ideal $I$ of a ring $R$ is a \sl
multiplicative-canonical \rm (or, simply, \sl m--canonical\rm) \sl
ideal of $R$ \rm if each regular fractional ideal $J$ of $R$ is
$I$--reflexive, i.e.  $J = (I:(I:J)) \cong
\mbox{Hom}_R(\mbox{Hom}_R(J, I), I)$.

It turns out that the previous question has a positive answer if
we assume a stronger condition on $I$: for each $n \geq1$, every
regular $R$--submodule of $R^n$ is $I$--reflexive. Under this
hypothesis we obtain that, for each $m \geq1$, every regular
$\RI$--submodule of $(\RI)^m$ is $\Hom_R(\RI,I)$--reflexive (see
Proposition \ref{Claim1-2} and Corollary \ref{canonical1}).
Moreover, $\Hom_R(\RI,I)$ is isomorphic to $\RI$ as an
$\RI$--module (see Theorem \ref{canonical2}). In particular, every
regular fractional ideal of $R\!\Join\! I$ is divisorial (see
Corollary \ref{cor:4.5}).

As a by-product, we obtain that, if $R$ is a Noetherian local
integral domain with an m--canonical ideal $I$,  then $R\!\Join\!
I$ is a reduced Noetherian local ring such that every regular
fractional ideal is divisorial (see Corollary \ref{cor:3.9}).




\section{Background on $\RI$}

Let $R$  be a commutative ring with unity, $T(R) \ (:=
\{\mbox{regular elements}\}^{-1}R)$ its total ring of fractions.
In this section we will give the definition of the ring $\RE$,
where $E$ is an $R$-submodule of $T(R)$ such that $E \!\cdot\! E
\subseteq E$ (note that this condition is equivalent to requiring
that there exists a subring $S$ of $T(R)$ containing $R$ and $E$,
such that $E$ is an ideal of $S$) and we will summarize some of
its properties we will need in this paper. For the sake of
simplicity, we will state these properties, for $E$ being a
nonzero (integral) ideal of $R$. \sl Mutatis mutandis \rm the
results hold in the general situation (cf.  \cite{d'a-f-1}, where
the interested reader can also find the details of the proofs).
\medskip

Let $\RE$ be the following subring of $R\times T(R)$ (endowed with
the usual componentwise operations):
$$
\RE:=\{(r,r+e) \ | \ r\ \in \ R, \ e  \in  E \} \ .
$$
It is obvious that, if in the $R$-module direct sum $R\oplus E$ we
introduce a multiplicative structure by setting $ (r,e)(s,f)
 := (rs,rf+se+ef)$, where $r,s \in R$
and $e,f \in E$, then we get the ring isomorphism $\RE \cong
R\oplus E$.


If $E=I$ is an ideal in $R$ (that we will assume to be proper and
different from $(0)$, to avoid the trivial cases), then the ring
$\RI$ is a subring of $R\times R$ and it is not difficult to see
that both the diagonal embedding $R \hookrightarrow \RI$ and the
inclusion $\RI \subset R\times R$ are integral. Moreover there
exist two distinguished ideals in $\RI$, $\boldsymbol{\mathfrak
O_1}  := (0) \times I$ and $\boldsymbol{\mathfrak O_2} := I \times
(0)$,  such that  $R \cong \RI/\boldsymbol{\mathfrak O_i}$, for
$i=1,2$.

As consequences of the previous facts we have:

\begin{prop}\label{bprop1} Let $I$ be a nonzero ideal of a ring $R$.

\begin{enumerate}

\item[\bf (1) \rm] If $R$ is a domain then $\RI$ is reduced and $
\boldsymbol{\mathfrak O_1} $ and $ \boldsymbol{\mathfrak O_2} $
are the only minimal primes of $\RI$.
\item[\bf (2) \rm ]   $R$ is reduced if and only if $\RI$ is reduced.
\item[\bf (3) \rm ]   $\dim(\RI) = \dim(R)$.
\item[\bf (4) \rm]   $R$ is Noetherian if and only if $\RI$ is Noetherian.
\end{enumerate}
 \end{prop}

Moreover it is possible to describe explicitly the prime spectrum
of $\RI$.

\begin{prop}\label{bprop2}

Let $P$ be a prime ideal of $R$ and set:
$$\aligned
\boldsymbol{\cal P}&:=  \{(p,p+i) \ | \ p \in P, \ i\in I \cap P \}, \\
\boldsymbol{\cal P}_{\!1}&:=\{(p,p+i) \ | \ p \in P, \ i\in I\},  \: \;  \text{and} \\
\boldsymbol{\cal P}_{\!2}&:=\{(p+i,p) \ | \ p \in P, \ i\in I\}.
\endaligned
$$\rm
\begin{description}
  \item[(a)] \it If $I \subseteq P$, then
  $\boldsymbol{\cal P}=\boldsymbol{\cal P}_{\!1}=\boldsymbol{\cal P}_{\!2}$
  is a prime ideal of $\RI$ and it is the unique prime ideal of $\RI$
  lying over $P$.\rm
  \item[(b)] \it If $I \nsubseteq
  P$, then $\boldsymbol{\cal P}_{\!1} \neq \boldsymbol{\cal P}_{\!2}$,
  $\boldsymbol{\cal P}_{\!1} \cap \boldsymbol{\cal P}_{\!2}=\boldsymbol{\cal P}$ and
  $\boldsymbol{\cal P}_{\!1}$ and $\boldsymbol{\cal P}_{\!2}$ are
  the only prime ideals of $\RI$ lying over $P$.\rm
  \item [(c)] \it The extension  $P(\RI)$  of $P$ in $\RI$  coincides with
  $\{(p,p+i) \ | \ p \in P, \ i\in IP\}$ and, moreover,
  $\sqrt{P(\RI)}=\boldsymbol{\cal P}$.
\end{description}
\it Furthermore, in case \rm (a) \it  we have:
$$R/P \cong
(\RI)/\boldsymbol{\cal P}  \mbox{ \;   and \;   } (\RI)_{
\boldsymbol{\cal P}} \cong R_P\!\Join\!I_P\,;$$ in case \rm (b)
\it  we have:
$$ R/P \cong (\RI)/\boldsymbol{\cal P}_{\!i} \mbox{  \:  and  \;  }
R_P \cong (\RI)_{\boldsymbol{\cal P}_{\!i }} \,, \mbox{  for
$i=1,2$}.$$ In particular, $R$ is a local ring if and only if
$\RI$ is a local ring.

\end{prop}

\section{Remarks on reflexivity}

We start this section recalling some definitions and results
related to the notion of multiplicative canonical ideal of a
domain, and giving the suitable generalizations for the
non--domain case.

Given an $R$--module $H$, for each $R$--module $F$, we can
consider the $R$--module $F^{\ast_H}:=\Hom_R(F, H). $  We have the
following  canonical homomorphism:
$$
\rho_F:  F \rightarrow  (F^{\ast_H})^{\ast_H}\,, \;\;\;\; a
\mapsto \rho_F(a)\,, \mbox{where  $\rho_F(a)(f) := f(a)$\,,}
$$
for all $f \in F^{\ast_H}$,\  $a \in F$.  We say
that the $R$--module $F$ is \it $H$--reflexive \rm (respectively,
\it $H$--torsionless\rm) if $\rho_F$ is an isomorphism
(respectively, monomorphism) of $R$--modules.\rm

Let $F$ be a regular $R$--submodule of $T(R)$ (i.e. $F$ contains a
$T(R)$--unit). It is not hard to prove that  each
$R$--homomorphism $F \rightarrow T(R)$ can be canonically extended
to an $R$--homomorphism $T(R) \rightarrow T(R)$.  Since
$\Hom_R(T(R), T(R))$ is canonically isomorphic to $T(R)$,  we have
that each $R$--homomorphism from $F$ into $T(R)$ is achieved by a
multiplication on $F$ by a unique element of $T(R)$.

 Given  a regular ideal $I$ of the ring $R$   and a
$R$--submodule $F$ of $T(R)$, set  $ (I:F) := \{z \in T(R) \mid zF
\subseteq I\}
 \cong \Hom_R(F,I)$.
If $F=J$ is a regular fractional ideal of $R$ then $(I:J)$ is
also a regular fractional  ideal of $R$.  Therefore, by the
previous considerations, we have a canonical isomorphism
  $(I :(I:J))  \xrightarrow{\sim}  (I:J)^{\ast_I}
\xrightarrow{\sim}  (J^{\ast_I})^{\ast_I}$. In this situation, we
can identify the map $\rho_J : J \rightarrow
(J^{\ast_I})^{\ast_I}$ with the inclusion $J \subseteq (I:(I:J))$,
so $J$ is $I$--torsionless.

We say that a regular ideal $I$ of a ring $R$ is \it  a
multiplicative-canonical ideal of $R$ \rm (or simply a \it
m--canonical ideal) \rm  if each regular fractional ideal $J$ of
$R$ is $I$--reflexive, i.e. the map $\rho_J : J \rightarrow
(J^{\ast_I})^{\ast_I}$ is an isomorphism or, equivalently, $J =
(I:(I:J))$.

Note that this definition is a natural extension of the concept
introduced in the integral domain case by W. Heinzer, J. Huckaba
and I. Papick \cite{hhp} and of  the notion of canonical ideal
given by J. Herzog and E. Kunz \cite[Definition 2.4]{hk} and by E.
Matlis \cite[Chapter XV]{matlis:1973} for 1--dimensional
Cohen-Macaulay rings. \ In general, given a Cohen-Macaulay local
ring $(R, M, k)$ of dimension $d$,  \it a canonical module \rm of
$R$  is an $R$--module $\boldsymbol{\omega}$ such that   the
$k$--dimension of $\mbox{Ext}^i_R(k, \boldsymbol{\omega})$ is 1
for $i=d$ and 0 for $i \neq d$.
If $R$ is not local, \it a canonical module for $R$ \rm  is an
$R$--module $\boldsymbol{\omega}$ such that all the localizations
$\boldsymbol{\omega}_{M}$ at the maximal ideals $M$ of $R$ are
canonical modules of $R_{M}$. When  a canonical module
$\boldsymbol{\omega}$ exists and it is isomorphic to an ideal $I$
of $R$, $I$ is called \it a canonical ideal  of $R$. \rm
In \cite[Proposition 4.3]{hhp} it is shown that a Noetherian
domain with dimension bigger than $1$ does not admit a
m--canonical ideal, while there exist (Noetherian) Cohen-Macaulay
domains of dimension bigger than $1$ with canonical ideal (e.g. a
Noetherian factorial domain $D$ of dimension $\geq 2$; in this
case, $D$ is a Gorenstein domain \cite[ Corollary 3.3.19]{b-h}).
Hence in higher dimension the notions of canonical ideal and
m--canonical ideal do not coincide.

\medskip  The following proposition extends outside of the
integral domain setting some results proved in \cite[Lemma 2.2
(a), (c) and Proposition 5.1]{hhp}.  The proof is standard and we
omit the details.

\begin{prop} \label{basic_can} Let  $I$ be a m--canonical ideal of a
ring $R$; then we have:

\begin{description}
\item[(a)]  $(I:I)= R  \cong \Hom_R(I, I)$ (the isomorphism
is realized by the canonical multiplication map $R \rightarrow \Hom_R(I, I)$).

\item[(b)]  If $L$ is an invertible ideal of $R$ (i.e. a regular
ideal such that $LL^{-1}=R$, where $L^{-1} := (R:L)$), then $IL$
is also a m--canonical ideal of $R$; in particular, for each
regular element $x \in R$, the ideal $xI$ is also a m--canonical
ideal of $R$.

\item[(c)]  Let $S$ be an overring of $R$, $R \subseteq S \subseteq
T(R)$, such that $(R:S)$ is a regular ideal of $R$. Then $(I:S) \
(\cong \Hom_R(S, I) =S^{\ast_I})$ is a m--canonical ideal of $S$.
\end{description}
\end{prop}

%
%
\medskip

We recall that a regular fractional ideal $J$ of a ring $R$ is
called a \it divisorial \rm ideal of $R$  if $(R:(R:J)) =J$.
Clearly, an invertible ideal of $R$ is a divisorial ideal. If
every regular fractional ideal of $R$ is divisorial, then $R$
itself is an m--canonical ideal.

The goal of the remaining part of this paper  is to study when
every regular fractional ideal of $\RI$ is divisorial.
We start by studying  some reflexivity properties related to the
notion of m--canonical ideal, in order to find an $\RI$--module
$E$, with the property that every regular ideal of $\RI$ is
$E$--reflexive.

\medskip
 Let $I$ be a regular ideal of a ring $R$ and
set
$$\aligned
 \calF_1 := \calF_1(R)   :=& \{ F \mid  F \mbox{ is a regular
 ($I$--torsionless)
$R$--submodule of $R$}  \} =\\
=& \{J \mid  J \mbox{ is a regular ideal of $R$ }  \} \ ;
\endaligned
$$
then we say that \it  the ring  $R$ \rm is  \it $(I,
\calF_1(R))$--reflexive \rm  (or, simply, \it $(I,
\calF_1)$--reflexive\rm) if each $F$ in $\calF_1(R)$ is
$I$--reflexive.
It is obvious that $R$ is  $(I, \calF_1)$--reflexive if and only
if $I$ is  an m-canonical ideal of $R$. (Note that  each regular
fractional ideal $J$ is $I$--reflexive if and only if $dJ$ is
$I$--reflexive, for each regular element $d \in R$  such that  $dJ
\subseteq R$.)

\medskip
Let $I$ be a regular ideal of a ring $R$. We have already observed
that every regular ideal of $R$ is $I$--torsionless. This property
holds more generally for every regular $R$--submodule of $R^n$,
for each $n \geq 1$. In other words, if $F$ is a regular
$R$--submodule of $R^n$ and if $x \in F \setminus \{0\}$, then we
can find  $h \in \Hom_R(F, I)$ such that $h(x) \neq 0$. As a
matter of fact,  write $x = (x_1, x_2, \dots, x_n)$ with $x_i\neq
0$, and let $\pi_i: R^n \rightarrow R$ be the projection on the
$i$-th coordinate. Choose $y \in I$ such that $x_iy \neq 0$, and
take $h\in \Hom_R(F, I)$ to be the composition $F\subseteq R^n
\xrightarrow{\pi_i} R\xrightarrow{y} I $. Set:
 $$
 \begin{array}{rl}
 \calF :=  \calF(R) := \{ F \mid & F \mbox{ is a regular,
 ($I$--torsionless) $R$--submodule of $R^n$,} \\
  & \mbox{  for some $n \geq1$} \}\,.
  \end{array}
 $$
 We say that \it the ring $R$ is $(I, \calF(R))$--reflexive \rm
 (or, simply, \it $(I, \calF)$--reflexive\rm)
if every  $F \in  \calF(R)$  is $I$--reflexive   (i.e. the
canonical monomorphism $\rho_F:  F \rightarrow  \Hom_R(\Hom_R(F,
I), I)$  is an isomorphism of $R$--modules).

 Note that if $R$ is  $(I, \calF)$--reflexive   then $I$ is a
m--canonical ideal of $R$, since each regular ideal $J$ of $R$
belongs to $ \calF(R)$.

\begin{prop} \label{Claim1-2}
Let $R$ be a ring admitting  a regular ideal  $I$ such that $R$ is
$(I,\calF(R))$--reflexive,   let $T$ be a subring of $R^m$, for
some $ m\geq 2$, containing the image of $R$ under the diagonal
embedding and set $I^{T} := \Hom_R(T, I)$. Let $E$ be any
$T$--module. Then the following  canonical maps are isomorphisms
of $T$--modules:
$$
\aligned \Hom_T(E, I^{T}) &\cong  \Hom_R(E, I)\\
\Hom_T( \Hom_T(E, I^{T}), I^{T}) &\cong \Hom_R(\Hom_R(E, I),
I))\,.
\endaligned
$$
\end{prop}

\noindent{\bf Proof.} Let  $E$ be a $T$--module. We can consider
$E$ as an $R$--module and so, by the ``Hom-tensor adjointness'',
we have that the map $\Hom_R(E, I^{T}) \rightarrow
\Hom_R(E\otimes_RT, I) = \Hom_R(E, I)$ defined by $ h \mapsto h',\
\mbox{where } h'(e):= h(e)(1 )$, for all $e \in E$, establishes an
isomorphism of $R$--modules.   On the other hand, note that $I^T$
is endowed with a structure of $T$--module, by setting  $z\cdot
f(t):= f(zt)$, for each $f \in I^{T}$ and $z, t \in T$; similarly
$\Hom_R(E, I)$ is a $T$--module, by setting  $z\cdot h' (e) := h'
(ze)$, for all $e\in E,\ z \in T,\ h' \in \Hom_R(E, I)$.

From the previous remarks it follows easily that the map:
$$\Phi: \Hom_T(E, I^{T}) \rightarrow \Hom_R(E, I)\,,\;\; h \mapsto
h'\,,\;\; \mbox{ where } h'(e):= h(e)(1 )\,,\; e \in E\,$$ is
bijective and preserves the sums. Moreover,  $\Phi$ is
$T$--linear, since $ \Phi(zh)(e) = (zh(e))(1)   = h(ze)(1) =
\Phi(h)(ze)   = ( z\Phi(h))(e)$, for all $e \in E,\ z \in T,\ h
\in  \Hom_T(E, I^{T})$.  Therefore  the map  $\Phi$ establishes an
isomorphism of $T$--modules.



By the  previous  isomorphism it follows that  the canonical maps
establish the following isomorphisms (as $T$--  and $R$--modules):
  $$\Hom_T( \Hom_T(E, I^{T}), I^{T}) \cong \Hom_R(\Hom_T(E,
I^{T}), I) \cong\Hom_R(\Hom_R(E, I), I))\,.$$ \hfill $\Box$

 \begin{cor} \label{canonical1}
 Let $R$ be a ring admitting  a regular ideal  $I$ such that
$R$ is  $(I,\calF(R))$--reflexive,   let $T$ be a subring of
$R^m$, for some $ m\geq 2$, containing the image of $R$ under the
diagonal embedding and set $I^{T} := \Hom_R(T, I)$.
Then every regular $T$--submodule of $T^n$, for some $n \geq1$, is
$I^T$--reflexive.
\end{cor}

\noindent{\bf Proof.} From the previous proposition, it follows
that  $E$ is a $T$--module $I^{T}$--torsionless (respectively,
$I^{T}$--reflexive) if and only if  $E$ is $I$--torsionless
(respectively, $I$--reflexive) as a $R$--module. Moreover, if $E$
is a regular $T$--submodule of $T^n$, then clearly $E$ is regular
$R$--submodule of $R^{mn}$. The conclusion is now straightforward.
%
\hfill $\Box$

\medskip
Notice that, in general, $I^T$ is not isomorphic to an ideal of
$T$.  However,  we will see in the next section that $I^T$ is
isomorphic to $T$ when  $T=\RI$ and $I$ is m--canonical.
 \section{$R \Join I$ when $I$ is a m--canonical ideal}

In this section we will investigate the construction $\RI$ in case
$I$ is an m--canonical ideal.  In particular we will extend,  to
not necessarily Noetherian rings,  one of the main results
obtained in \cite[Theorem 11]{d'a}.

\begin{thm}\label{canonical2}
Let $I$ be an ideal of a ring $R$ such that the canonical
(multiplication) map $R \rightarrow \Hom_R(I,I)$ is an isomorphism
(e.g. let $I$ be a m--canonical ideal of $R$; Proposition
\ref{basic_can} (a)).   Then   $ \Hom_R(\RI, I)$ is isomorphic as
$\RI$-module to $\RI$.
\end{thm}

\noindent {\bf Proof.} 
Since $I \cong \Hom_R(R, I)$ (under the map $\iota  \mapsto \iota
\cdot \mbox{--}$, for $\iota \in I$) and $R \cong \Hom_R(I, I)$
(under the map $x \mapsto x\cdot \mbox{--}$, for  $x \in R$), we
deduce immediately that there is a canonical isomorphism of
$R$--modules $\RI \cong R \oplus I \cong \Hom_R(R, I) \oplus
\Hom_R(I, I)$. Moreover we have the  following  canonical
isomorphism of $R$--modules:
$$
\Hom_R(R, I) \oplus  \Hom_R(I,I) \rightarrow  \Hom_R(R\!\Join \!
I, I)\,,\;\;  (g_1, g_2) \mapsto g\,,
$$
(where $g: R\!\Join \!I \rightarrow I, \ (z, z + j) \mapsto g_1(
z)+g_2(j)$, for each $z \in R$ and $i\in I$).

Note that the composition map  $$\RI \rightarrow \Hom_R(R\!\Join
\! I, I)\,,\;\;  (x,x+i) \mapsto g_{(x,i)},$$ where
$g_{(x,i)}((z,z+j)) := iz+xj$ (for all $ x, z \in R$ and $i,j \in
I$) is obviously  an $R$--isomorphism, but it is not an
$\RI$--isomorphism.

In order to  get an $\RI$--isomorphism we  consider the following
map:
$$\sigma: \RI \rightarrow
\Hom_R(R\!\Join \!I, I)\,,\;\;(x,x+i) \mapsto f_{(x,i)}$$ where
$f_{(x,i)}(z, z+j) :=  xj + (z+j)i $, for all $ x, z \in R$ and
$i,j \in I$.

It is not difficult to check that $\sigma$ is an injective
$\RI$--homomorphism (recall that the natural structure of
$\RI$--module on $\Hom_R(\RI, I)$, is defined by the scalar
multiplication by $ (x', x' +i') \cdot f((z, z+j)):= f((x',
x'+i')(z, z+j))$).

It remains to prove that $\sigma$ is surjective, that is: for each
$f \in   \Hom_R(R\!\Join \!I, I)$, there exists $(\bar{x},\bar{x}+
\bar{\iota}) \in \RI$, such that $f = \sigma((\bar{x},\bar{x}+
\bar{\iota})) = f_{(\bar{x}, \bar{\iota})}$\,. Let $\bar{\iota} :=
f((1,1))$ and let $\bar{x} := (f((-\bar{j }, 0))/\bar{j }) $, for
some  nonzero (regular) element $\bar{j }\in I$.  Note that
$(f((-\bar{j }, 0))/\bar{j }) $ does not depend on the choice of
$\bar{j }$, since $j' (f((-j, 0))= f((-jj', 0)) = j (f((-j', 0))$,
i.e. $(f((-j, 0))/j) =(f((-j', 0))/j') $, for any two nonzero
(regular) elements $j, j' \in I$. The previous relation shows also
that $(f((-\bar{j }, 0))/\bar{j })I \subseteq I$,   therefore from
the canonical isomorphism of $R$-modules $\Hom_R(I, I) \cong R$,
we deduce that $(f((-\bar{j }, 0))/\bar{j })$ (which a priori is
an element of $T(R)$) belongs to $ R$. Note also that, for each $j
\in I$, $ f((0, j))= f((-j, 0)) +jf((1,1)) = (\bar{x }  +
\bar{\iota})j$. Therefore, for each $z \in R$ and $j \in I$, we
have:
$$\begin{array}{rl}
f((z, z+j)) &=f((z,z)) + f((0, j)) =\\
&= zf((1,1))+ jf((1,1)) +  f((-j, 0)) = \\
&= (z+ j)\bar{\iota } + \bar{x } j =\\
&= f_{(\bar{x } , \bar{\iota})} ((z, z+j)) \,.
\end{array}
$$

Hence we can conclude that the map $\sigma: \RI \rightarrow $
$\Hom_R(R\!\Join \! I, I) $ is an isomorphism of $\RI$-modules.
 \hfill $\Box$

\medskip

We remark that an alternative proof of the previous result can be
given by showing that $\Hom_R(\RI,I)$  is a free $\RI$--module of
rank one: if we denote by $\pi: \RI \rightarrow I$ the canonical
projection, $(r, r+ i) \mapsto i \ (= (r+i) - r)$, then it is
possible to show that $\{\pi\}$ is a basis for $\Hom_R(\RI,I)$ as
a $\RI$--module.

\begin{cor} \label{cor:4.5} Let $R$ be a ring admitting a
regular ideal $I$ such that $R$ is  $(I, \calF(R))$--reflexive.
Then $\RI$ is  $(\RI, \calF(\RI))$--reflexive. In particular,
every regular fractional ideal of $\RI$ is divisorial.
\end{cor}
\noindent{\bf Proof.} If we set $T:=\RI$, by Theorem
\ref{canonical2} we have $I^{T} \cong T$; moreover by  Proposition
\ref{canonical1},  $T$ is  $(I^{T}, \calF(T))$--reflexive.
\hfill $\Box$

\medskip

It is natural to ask whether the last statement of the previous
Corollary \ref{cor:4.5} holds by assuming that $I$ is  an
m--canonical ideal  if $R$.
A related problem  is to find conditions on
$R$ so that if $
I$ is a m--canonical ideal of $R$, then $R$ is
$(I, \calF)$--reflexive.   The remaining part of this section is
an investigation in this direction.

Recall that \it a Marot ring \rm is a ring such that each regular
ideal is generated by its set of regular elements and a \it ring
has few zero divisors \rm if the set of zero divisors is a finite
union of prime ideals \cite[page 31]{h}. Recall that a Noetherian
ring is always a ring with few zero divisors and a ring with few
zero divisors is a Marot ring; moreover an overring of a Marot
ring is a Marot ring \cite[Theorem 7.2 and Corollary 7.3]{h}.

The following result extends  \cite[Proposition 3.6]{hhp} to the
non integral domain case   and shows that the conclusion of the
previous Corollary \ref{cor:4.5} can be characterized in several
ways in the context of Marot rings.

\begin{prop} \label{can-div}  Given a Marot ring $R$, the following
statements are equi\-va\-lent:
\begin{description}

\item[(i)] $R$ has a principal m--canonical ideal;

\item[(ii)]  $R$ has an invertible m--canonical ideal;

\item[(iii)]  $R$ has a divisorial m--canonical ideal;

\item[(iv)]  each regular  fractional ideal of $R$ is divisorial.

\end{description}\end{prop}
 \noindent{\bf Proof.}  It is obvious that
(iv)$\Rightarrow$(i)$\Rightarrow$(ii)$\Rightarrow$(iii).

(iii)$\Rightarrow$(iv).   We start by showing a generalization of
\cite[Lemma 3.1]{hhp}:

\bf Claim 1. \it Given two regular fractional ideals $I, J$ of $R$,
then:
$$
(I:(I:J)) = \bigcap \{zI \mid z \mbox{ is a regular element of
$T(R)$ with } J \subseteq zI \}\,.
$$ \rm

Let $x \in (I:(I:J))$ and let $z \in T(R)$ be a regular element
such that  $ J \subseteq zI$. Then, clearly, $x (I:J) \subseteq I$
and $z^{-1}J \subseteq I$. Therefore  $xz^{-1} \in I$, i.e. $x \in
zI$, and so $(I:(I:J))\subseteq  \bigcap \{zI \mid z \mbox{ is a
regular element of  $T(R)$ with } J \subseteq zI \}$.  On the
other side, let $x \in  \bigcap \{zI \mid z \mbox{ is a regular
element of  $T(R)$ with } J $ $\subseteq zI \}$.
If $u$ is  a
regular element of $(I:J)$,  then  $J \subseteq
u^{-1}I $. Therefore, $x \in  u^{-1}I $, i.e. $xu \in I$.
Since $(I:J)$ is
a regular fractional ideal of the Marot ring $R$, it follows that
$x \in (I:(I:J))$.

Since $(I:J)$ is
a regular fractional ideal of the Marot ring $R$, if $u$ is  a
regular generator of $(I:J)$,  then we have that $J \subseteq
u^{-1}I $. Therefore, $x \in  u^{-1}I $, i.e. $xu \in I$, and so
$\bigcap \{zI \mid z \mbox{ is a regular element of  $T(R)$ with }
J \subseteq zI \} \subseteq (I:(I:J))$.

\bf Claim 2. \it  If $\{ J_{\alpha} \mid \alpha \in A \}$ is a
family of divisorial ideals of $R$ such that  $\bigcap_\alpha
J_\alpha\neq
  (0)$, then $\bigcap_\alpha J_\alpha$ is divisorial. \rm

Let $\{z_\beta\mid \beta \in B\}$ be a family of regular elements
in $T(R)$ such that $\bigcap _ \beta z_ \beta R \neq   (0)$. Note
that, for each $z_ \beta$,  we have $(R:z_ \beta R)  \subseteq (R
:(\bigcap _ \beta z_ \beta R))$ and so $\sum _ \beta (R:z_ \beta
R)  \subseteq (R :(\bigcap _\beta z_ \beta R))$.  Moreover, it is
easy to see that, given a family of regular fractional ideals
$\{L_\beta \mid \beta \in B\}$ of $R$:
$$ (R : \sum_\beta L_\beta) =  \bigcap_\beta (R: L_\beta)\,.
$$
Therefore:
$$
\begin{array}{rl}
\bigcap _ \beta z_ \beta R \subseteq&  (R:(R:(\bigcap _ \beta z_
\beta R))) \subseteq
(R:\sum _ \beta (R:z_ \beta R)) =\\ =&  \bigcap _ \beta  (R: (
(R:z_\beta R))) =  \bigcap_ \beta z_ \beta R\,.
\end{array}$$
 Hence $\bigcap _ \beta z_ \beta R $ is a divisorial ideal of
$R$.

The conclusion follows easily, since by Claim 1 (for $I=R$), a
divisorial ideal in a Marot ring  is the intersection of a family of
principal regular fractional ideals.

\bf Claim 3. \it Given two regular fractional ideals $I, J$ of $R$,
assume that $I$ is divisorial. Then
$
(I:(I:J))$ is divisorial. \rm

This is an easy consequence of Claims 1 and 2,  since $I$ divisorial
implies that $zI$ is divisorial, for each regular element  $z \in
T(R)$.

 Now we can easily conclude the proof since, if $J$ is a regular
fractional ideal of $R$ and $I$ is the divisorial m--canonical
ideal of $R$, then $J = (I:(I:J))$ and, by Claim 3,  $ (I:(I:J))$
is divisorial. \hfill $\Box$

\begin{rem} \rm  \bf   (a)  \rm Note that the hypothesis that $R$ is
a Marot ring is essential in Claim 1
of the
proof of Proposition \ref{can-div}, see \cite[Theorem 8.3 and Section
27, Example 11]{h}. Several classes of examples of Marot rings are
given in  \cite[Section 7]{h}.

\noindent \bf (b)  \rm Note that if $R$ is a ring with  few
divisors and $I$ is an  ideal of $R$ then $\RI$  has few zero
divisors and so is a Marot ring.

\end{rem}

Let $R$ be  an integral domain and $I$ a nonzero  ideal of $R$.
Recall that  $R$ is said to be \it $I$--reflexive \rm
(respectively, \it $I$--divisorial\rm) in the sense of Bazzoni and
Salce \cite{bs} (cf. also \cite{O}), if  each $I$--torsionless  $\Hom_R(I,I)$--module
of finite rank (respectively, of rank 1)  is $I$--reflexive.




\begin{prop} \label {prop:3.18}
Let  $R$ be  an integral domain and $I$ a nonzero  ideal
of $R$.
Then
$R$ is  $(I, \calF)$--reflexive  (respectively,   $(I, \calF_1)$--reflexive)
if and only if $R$ is $I$--reflexive (respectively, $I$--divisorial) and
$R = (I:I)$.
\end{prop}

\noindent{\bf Proof.} First note that, if  $I$ is  a m--canonical
ideal of $R$ (i.e. if  $(I, \calF_1)$--reflexive; this happens
when $R$ is  $(I, \calF)$--reflexive),  then $\Hom_R(I,I) \cong R$
(Proposition \ref{basic_can} (a)).

If  $R$ is  $(I, \calF)$--reflexive (respectively,
$(I, \calF_1)$--reflexive)  we need to verify that each
$I$--torsionless  $R$--module $G$  of finite rank (respectively,
of rank 1)  is $I$--reflexive.  By \cite[Lemma 5.1]{f-s}, such  a
$G$ can be embedded in $I^n$, where $n$ is the rank of $G$, hence
$G$ belongs to  $\calF$ (respectively, $\calF_1$)  and so
$G$ is     $I$--reflexive.

Conversely, let $R = (I:I)$ and let $F \in  \calF $ (respectively, $F \in
 \calF_1$).  For  each nonzero element $i \in I$, then
$iF$  is a  $I$--torsionless  $R$--module  of finite rank
(respectively, of  rank 1) by \cite[Lemma 5.1]{f-s}, since $iF
\subseteq iR^n \subseteq I^n$ (respectively,  $iF \subseteq iR
\subseteq I$).  Therefore   $iF$ is $I$--reflexive  and so  $F$ is
$I$--reflexive. \hfill $\Box$
\medskip

\begin{cor} \label{cor:3.9}
Let $R$ be a Noetherian local integral domain and let $I$ be a
m--canonical ideal of $R$ and  set $T:= R\!\Join\! I$.
Then $T$ is a Noetherian local reduced ring, with $\dim(T) = \dim(R)$, such that every
regular fractional ideal of $T$ is divisorial.
\end{cor}

\noindent{\bf Proof.} By Propositions \ref{bprop1} and
\ref{bprop2} we know that $T$ is a Noetherian local reduced ring
with $\dim(T) = \dim(R)$.  Note that $I$ is a m--canonical ideal
of an integral domain $R$ if and only if
$I$--divisorial  and $(I:I)=R$ (Proposition \ref{prop:3.18}).
Moreover, in the   Noetherian local integral domain case, if
$(I:I)=R$,  then  $R$ is $I$--divisorial if and only if $R$ is
$I$--reflexive, by Bazzoni's generalization of Matlis' 1-dimensional
theorem \cite[Theorem 3.2]{b}.   By reapplying Proposition
\ref{prop:3.18}  we know, in this case,  that $R$ is
$I$--reflexive and $(I:I)=R$ if and only if $R$ is
 $(I, \calF(R))$--reflexive.   The conclusion follows immediately from
Corollary \ref{cor:4.5}.  \hfill $\Box$

\smallskip
Note that the assumption that a Noetherian  domain $R$ admits  an
m-canonical ideal implies dim$(R) \leq 1$ (by \cite[Proposition
4.3]{hhp}). Therefore, under the assuntions of Corollary
\ref{cor:3.9}, we can conclude that $\RI$  is a one-dimensional
reduced Gorenstein local ring \cite[Korollar 3.4]{hk}.

\bigskip
\noindent \bf Acknowledgments. \rm
The authors wish to express their sincere gratitude to the referee  for his/her valuable comments and helpful suggestions.

\noindent
Marco D'Anna\\ Dipartimento di Matematica e Informatica\\
Universit\`a di Catania\\ Viale Andrea Doria 6\\ 95125 Catania,
Italy\\
{\texttt {mdanna@dmi.unict.it}} \\

\smallskip\noindent
Marco Fontana\\ Dipartimento di Matematica\\
Universit\`a di Roma degli Studi ``Roma Tre''\\ Largo San Leonardo
Murialdo 1\\ 00146 Roma, Italy\\
\texttt{fontana@mat.uniroma3.it}

\end{document}